\documentclass[10pt]{amsart}

\usepackage{graphicx}
\usepackage[english]{babel}
\usepackage{latexsym}
\usepackage{amsmath}
\usepackage{amssymb}
\usepackage{mathrsfs}
\usepackage{graphicx}
\usepackage{xcolor}
\usepackage{bbm}
\usepackage{bbold}
\input epsf

\newcommand{\cl}{\operatorname{cl}}
\newcommand{\co}{\operatorname{co}}
\newcommand{\cone}{\operatorname{cone}}
\newcommand{\core}{\operatorname{core}}
\newcommand{\Cop}{\operatorname{Cop}}
\newcommand{\dom}{\operatorname{dom}}
\newcommand{\Fd}{\operatorname{Fd}}

\newcommand{\mix}{\operatorname{mix}}
\newcommand{\op}{\operatorname{op}}
\newcommand{\olim}{\mathop{o\text{-}\lim}}
\newcommand{\Orth}{\operatorname{Orth}}
\newcommand{\rlim}{\mathop{r\text{-}\lim}}
\newcommand{\subd}{\underline{\partial}}
\newcommand{\supd}{\overline{\partial}}
\newcommand{\tsubd}{\underline{\partial}{}^c}
\newcommand{\tsupd}{\overline{\partial}{}^c}
\newcommand{\Sbl}{\mathop{\fam0 Sbl}}
\newcommand{\CSc}{\mathop{\fam0 CS}\nolimits_c}

\newcommand{\QL}{\mathop{\fam0 QL}}
\newcommand{\tCSc}{\mathop{\fam0 CS}\nolimits_c^c}
\newcommand{\tQL}{\mathop{\fam0 QL}\nolimits^c}
\newcommand{\Proclaim}[1]{\smallskip{\bf#1}~\sl}
\newcommand{\proclaim}[1]{{\bf#1}~\sl}
\newcommand{\Endproc}{\rm}
\begin{document}

\title{Quasidifferentials in Kantorovich Spaces}

\author{Basaeva E.K. \and Kusraev A.G. \and Kutateladze S.S.}



\dedicatory{In memory of Vladimir F. Demyanov}

\begin{abstract}
This is an overview of the quasidifferential calculus for
the mappings that arrive at Kantorovich spaces. The necessary optimality
conditions are also derived for multiple criteria optimization problems
with quasidifferentiable data.
\end{abstract}

\maketitle

 \section*{Introduction}
 \label{intro}
A mapping is {\it quasidifferentiable\/} at an interior point of
the domain of definition provided that there is a directional derivative at
the point which can be presented as difference of two sublinear operators.
In this event the quasidifferential is introduces by a natural extension of the
Minkowski duality. This leads to a rather broad class of the  mappings
admitting linearization which contains convex and concave operators.

 The central problem of the calculus consists in providing some formulas
for the quasidifferential  of a composite mapping. This problem splits into
the three stages: (1)---search for the explicit representation of
the directional derivative of the mapping through the directional derivatives
of the terms of the composition; (2)---representation of
the so-obtained directional derivative as difference of sublinear operators;
and (3)---calculation of the quasidifferential of  the composite through
the quasidifferentials of the terms of the composition.
The first stage consists in calculating the relevant limits
and uses the tools of the classical analysis with due technical modifications.
The second stage is either obvious or involves some artificial tricks.
The third stage bases on the Minkowski duality in the version extended to
the broader class of quasilinear operators, i.e., differences of sublinear
operators.

Using this approach we can derive all
the main formulas of
quasidifferential calculus, i.e., the quasidifferentials of a sum,
a product, a fraction, a composition, a supremum, and an infimum
(cp.~\S\,3).

 Quasidifferential calculus enables us to derive
the necessary optimality conditions with quasidifferentiable  constraints
given as
along the lines of  subdifferential calculus; cp.~\S\S\,4,\,5.
More details are collected in \cite[Chapter~5]{SBD}.)

 \section{Kantorovich Spaces}
 \label{sec:1}
 In this section we briefly outline
the main definitions and facts of the theory of Kantorovich spaces
which we will use later.
More details are collected in~\cite{15}, \cite{KaA}, \cite{45}, and \cite{DOP}.

 \paragraph{1.1.}
 In what follows, $\mathbb R$ is~the field of the reals.
An~{\it ordered vector space\/} (over $\mathbb R$) is some pair
 \mbox{$(E,\leq)$}, with $E$  a~real vector space and $\leq$
an~order on~$E$, satisfying the conditions:

 \textbf{(1)} if $x\leq y$ and $u\leq v$, then $x+u \leq y+v$ for all
 $x,y,u,v\in E$;

 \textbf{(2)} if $x\leq y$ then  $\lambda x\leq\lambda y$ for all
 $x,y\in E$ and $0\leq\lambda \in \mathbb R$.

 Introduction such a relation on a vector space
amounts to distinguishing some set $E_+\subset E$,
 called the {\it positive cone\/} in~$E$  and enjoying the
properties::
 $$
 E_++ E_+\subset E_+, \qquad \lambda E_+ \subset E_+\ \ (0\leq
 \lambda \in \mathbb R).
 $$
 In this event the quasiorder  $\leq$  and the positive cone
$E_+$ are connected as follows: $x\leq y\,\Longleftrightarrow\, y-x\in E_+$ $(x,y\in
 E)$. The elements of $E_+$ are referred to as {\it
 positive}. The quasiorder  $\leq$ is an order if and only if
  $E_+ $ is an {\it ordering\/} or {\it salient\/} cone, i.e.,
$E_+ \cap -E_+ =\{0\}$.

 \paragraph{1.2.}~An order vector space whose order makes it into a lattice
is called a~{\it vector lattice}.
 To every two elements  $x$ and $y$ of a vector lattice $E$
 there correspond their {\it join\/} or {\it supremum\/}
$x\vee y\!:=\sup\{x,y\}$ and their {\it meet\/} or {\it infimum\/}
 $x\wedge y\!:=\inf\{x,y\}$. In particular, each~$x$ in a vector lattice~$E$
has the {\it positive part\/} $x^+ := x\vee 0$, the
 {\it negative part\/} $x^- :=(-x)^+ :=-x \wedge
 0$, and the {\it modulus\/}~$|x| := x\vee (-x)$.

 A vector lattice $E$ is {\it Dedekind complete\/} whenever
 each nonempty order bounded subset of~$E$ has a supremum and an infimum.
A Dedekind complete vector lattice is called a~{\it Kantorovich space\/}
or, briefly, a~{\it  $K$-space\/} in memory of Leonid Kantorovich who was the
first to distinguish this class of space in~\cite{Ka}.
Below $E$ is an arbitrary $K$-space unless specified otherwise.

 \paragraph{1.3.}~Elements $x,y\in E$ are called {\it disjoint\/},
which is written as $x\perp y$, provided that $|x|\wedge
 |y| =0$. The set
 $$
 M^\perp\!:=\{x\in E:\,(\forall y \in M)\,x\perp y\},
 $$
 with $ M \subset E$, is the {\it disjoint complement\/} of~$M$.

 A~{\it band\/} of~$E$ is a set of the form
$M^\perp$, with  $M \subset E$ and $M\ne\varnothing$. The inclusion ordered
set of all bands of~$E$ is the complete Boolean algebra
 ${\mathfrak B}(E)$ with the Boolean operations
 $\vee$, the meet; $\wedge$, the join; and $(\cdot)^*$, the Boolean complement
acting as follows:
 $$
 L\wedge K=L\cap K,\ \ L \vee K=(L \cup  K)^{\perp\perp},\ \
 L^\ast=L^\perp\quad({ L, K \in {\mathfrak B}(E)}).
 $$
 The algebra ${\mathfrak B}(E)$ is called the {\it base\/} of~$E$.

 \paragraph{1.4.}
 A $K$-space $E$ is called {\it extended\/} or {\it universally complete\/}
provided that each nonempty set of pairwise disjoint elements of ~$E$
has a supremum. Universally complete $K$-spaces are isomorphic if and only if
so are their bases. An arbitrary $K$-space admits an embedding
in a universally complete $K$-space with the same base. This completion is
unique to within isomorphism. We now list the main examples
of universally complete $K$-spaces:

 \textbf{(1)}~The space $L^0(\Omega, \Sigma, \mu)$ of cosets of
almost finite measurable functions on $\Omega$, where $(\Omega,\Sigma,
 \mu)$~ is a~measure space and $\mu$ is $\sigma$-finite (or, which is more general,
has the direct sum  property; cp.~\cite{KaA}, \cite{109}.
 The base of $L^0(\Omega, \Sigma , \mu )$ is isomorphic with
$\Sigma/ \mu^{-1}(0)$, the Boolean factor algebra  of measurable sets over
negligible sets.

 \textbf{(2)}~The space $ C_{\infty} (Q)$ of continuous functions
on an extremally disconnected compact space $Q$
ranging in the extended real axis and each taking the value
 $\pm \infty $ only on a~nowhere dense subset; cp.~\cite{15} and~\cite{DOP}.
 The base of this $K$-space is isomorphic to the
Boolean algebra of clopen subsets of~$Q$.

 \paragraph{1.5.}~The order of a vector lattice provides a few types of convergence
 of nets and sequences. A net $(x_\alpha)_{\alpha\in\mathrm{A}}$ in~$E$ is {\it order convergent\/} or
 {\it $o$-converges\/} to $x\in E$ provided that there is a decreasing net
 $(e_\beta )_{\beta\in {\mathrm{B}}}$ in~$E$ such that
 $\inf\{ e_\beta :\, \beta\in \mathrm{B}\}=0$ and to each $\beta\in
 \mathrm{B}$ there is $\alpha(\beta)\in \mathrm{A}$ satisfying
  $|x_\alpha  - x|\leq e_\beta $ for all
 $\alpha(\beta)\leq\alpha \in \mathrm{A}$. In this event $x$
is the {\it order limit\/} or
 {\it \hbox{$o$-}limit\/} of $(x_\alpha)$  and we write
 $x=\olim x_\alpha$ or $x_\alpha\overset{(o)}\longrightarrow x$.

 If we replace $(e_\beta)$ by a sequence $(\lambda_n e)_{n\in\mathbb N}$, where $0 \leq e\in
 E_+$ and $(\lambda _n)_{n\in \mathbb N}$ is a~numeric sequence satisfying
 $\lim_{n\to\infty} \lambda _n = 0$, then we say that
  $(x_\alpha )_{\alpha \in \mathrm{A}}$
  (relatively uniformly) converges to $x\in E$ with  {\it regulator~$e$}. The element $x$
 is the {\it  $r$-limit\/} of $(x_\alpha )$. We use the notations
  $x = \rlim_{\alpha \in \mathrm{A}} x_\alpha$
 and $x_\alpha \overset{(r)}\longrightarrow x$.

 \paragraph{1.6.}~Let $E$ and $F$ be~vector lattices. By
 $L(E,F)$ we will denote the space of linear operators from~$E$ to~$F$.
 We call a linear operator $T:E\to F$
 {\it positive} in case
 $T(E_+)\subset F_+$; {\it regular},
in case $T$ can be presented as difference of positive operators;
{\it order bounded\/} or, briefly, {\it $o$-bounded} in case
 $T$ maps each order bounded subset of~$E$ to an order
 bounded subset of~$F$.

Every positive operator is obviously order bounded; hence,  so
is the difference of order bounded operators. Therefore, every
regular operator is order bounded. The converse, failing in general,
is valid provided that $F$ is Dedekind complete, as follows from
the definitive Riesz--Kantorovich Theorem.

 \paragraph{1.7.} \proclaim{Riesz--Kantorovich Theorem.}
 Let $E$ be a~vector lattice and let $F$ be a~$K$-space.
 The set of all order bounded operators $L^\sim(E,F)$,
 ordered by the positive cone $L^\sim(E,F)_+$,
 is a~$K$-space. \Endproc

 The basics of the theory of regular operators in~$K$-spaces were
 laid by L.~V. Kantorovich in~\cite{Ka} where
 the above theorem had appeared firstly. F.~Riesz formulated a similar theorem for the space of continuous linear functionals on the vector lattice $C([a,b])$ in
  his celebrated talk \cite{Ris} at the International Mathematical Congress in Bologna in 1928.

 \paragraph{1.8.}
 A linear operator $T:E\to F$ is called
 {\it order continuous\/} (respectively, {\it sequentially
 $o$-continuous\/} or {\it order $\sigma$-continuous}), provided that $Tx_\alpha
 \overset{(o)}{\longrightarrow} 0$ in $F$ for every net $(x_\alpha)$
 $o$-converging to the zero of~$E$ (respectively, $Tx_n
 \overset{(o)}{\longrightarrow} 0$ in~$F$ for every sequence
  $(x_n)$ that $o$-converges to the zero of~$E$). The sets of all
  order continuous and order $\sigma$-continuous operators from~$E$ to~$F$
are denoted by $L_n^\sim(E,F)$  and
 $L_\sigma^\sim (E,F)$.

 \paragraph{1.9.}
 \proclaim{Theorem.}
 Let $E$ and $F$ be~vector lattices, with $F$ Dedekind complete.
Then $L_n^\sim(E,F)$  and $L_\sigma^\sim (E,F)$  are bands of~$L^\sim(E,F)$. \Endproc

 \paragraph{1.10.}~Consider some $K$-space~$E$. Let
 $L^\sim (E)\!:=L^\sim(E,E)$ be the~space of regular endomorphisms
 of~$E$. Denote by $\Orth(E)$ the least band of $L^\sim (E)$ which contains the identity operator
 $I_E$; i.e., $\Orth(E):=\{I_E\}^{\perp\perp}$.
 The elements of $\Orth(E)$ are called {\it orthomorphisms}.
 The orthomorphisms enjoy  a few remarkable properties:

 \textbf{(1)}~$T\in L^\sim (E)$ is an orthomorphism if and only if
 $T$ commutes with order projections; i.e., $\pi T=T\pi $ for all $\pi\in\mathfrak P(E)$.

 \textbf{(2)}~Every orthomorphism is order continuous.

 \textbf{(3)}~Every two orthomorphisms commute.

 \textbf{(4)}~If $\pi,\rho\in \Orth(E)$, then $\rho\circ \pi\in \Orth(E)_+$.

 \textbf{(5)}~The $K$-space $\Orth(E)$ with composition as multiplication
 is an~{\it $f$-algebra}; i.e., $\Orth(E)$ is a~vector lattice and a commutative algebra
 such that
 $$
 \alpha, \pi, \rho\in \Orth(E)_+,\ \, \pi\perp \rho \Longrightarrow
 (\alpha \pi) \perp \rho.
 $$

 In what follows $A:=\Orth(E)$ is the  orthomorphism algebra of~$E$.

 \paragraph{1.11.}
 Let $X$ be a~vector space, while $E$ is a~$K$-space. An operator $p:X\to E$ is {\it sublinear\/}
 provided that $p(x+y)\leq p(x)+p(y)$ and $p(\lambda x)=\lambda p(x)$ for all $x,y\in X$ and $0\leq\lambda\in\mathbb{R}$.
 Denote the set of all sublinear operators from $X$ to~$E$ by $\Sbl(X, E)$.
 Define addition and multiplication on $\Sbl(X,E)$ pointwise, i.e.,
 $(p+q)x\!:=p(x)+q(x)$ and  $p\leq q\,\Longleftrightarrow\,(\forall x\in X)\,p(x)\leq q(x)$ $(x\in X)$.
 Clearly, the addition and order agree with one another as usual: $p_1\leq p_2\,\Longrightarrow\,p_1+q\leq p_2+q$.

 To each  $p\in\Sbl(X,E)$ we can uniquely assign its {\it support set\/}
 $$
 \partial p\!:=\{T\in L(X,E):\ Tx\leq p(x)\ (x\in X)\}.
 $$
 Put $\CSc(X,E):=\{\partial p:\,p\in\Sbl(X,E)\}$.

 \paragraph{1.12.}
 Denote by $\QL(X,E)$ the set of all mappings from
 $X$ to $E$ that are representable as differences of sublinear operators.
 Under the natural algebraic operations, $\QL(X,E)$ becomes an $A$-module. Indeed,
 if  $f=p-q$ for some $p,q\in\Sbl(X,E)$ and
 $\alpha\in\Orth(E)$ then
 $$
 \alpha f:=\alpha\circ f=(\alpha^+p+\alpha^-q)-(\alpha^-p+\alpha^+q),
 $$
 and so $\alpha f\in\QL(X,E)$.
 The elements of $\QL(X,E)$ are called {\it quasilinear operators}.

 In other words, $\QL(X,E):=\Sbl (X,E)-\Sbl (X,E)$ and the structure of
an ordered $A$-module is induced from~$E^X$, i.e., this is done pointwise.
In particular, the order on $\QL(X,E)$ is determined from the positive cone
 $\{f\in \QL(X,E):\, f(x)\geqslant 0$ $(x\in X)\}$.

 To each pair of sublinear operators $p,q\in\Sbl (X,E)$ we put into correspondence
 the quasilinear operator $\phi(p,q):x\mapsto p(x)-q(x)$ $(x\in X)$.
 Observe that the pairs $(p,q)$ and $(p',q')$
 yield the same quasilinear operator, i.e. $\phi(p,g)=\phi(p',g')$, provided that
  $p-q=p'-q'$. Wi will call these pairs {\it equivalent}.
  The coset of $(p,q)$ will be denoted by $[p,q]$.

 \paragraph{1.13.}
 To each pair $p,q\in\Sbl(X,E)$ we can assign the pair of their support sets
 $\partial p$ and $\partial q$.
 Now, to each quasilinear operator  $l=p-q$ there corresponds
 the coset $[p,q]$ as well as the coset $[\partial p,\partial q]$.

 Finally, if $l=p-q$ then we can put $\mathcal D l:=[\partial p,\partial q]$
 and observe that $\mathcal D l$ does not depend on the particular representation
 of $l$ as difference of sublinear operators.
 The element $\mathcal D l$ of the $A$-module $\CSc (X,E)$ is called
 the  {\it quasidifferential\/} of $l$ (at the zero).
 Moreover, for the support sets $\partial p$ and $\partial q$ we use the following
 terms  and denotations: $\subd l:=\partial p$ is the {\it subdifferential\/} of
 ~$l$ (at the zero) and $\supd l:=\partial q$ is the {\it superdifferential\/} of~$l$
 (at the zero).

 \paragraph{1.13.}
 Let us inspect what happens to various operations under the mapping $\mathcal D : l\mapsto \mathcal D$.
 If  $\alpha\in\Orth(E)$ and $l,l_1,\dots,l_n\in\QL(X,E)$ then
 $$
 \aligned \mathcal D(\alpha l)=\alpha\mathcal Dl &= [\alpha^+\subd
 l+\alpha^-\supd l,\, \alpha^-\subd l+\alpha^+\supd l\/],
 \\
 \mathcal D(l_1+\ldots+l_n)&=\mathcal Dl_1+\ldots+\mathcal Dl_n=
 \\
 &=[\,\subd l_1+\ldots+\subd l_n,\,\supd l_1+\ldots+\supd l_n\/],
 \endaligned
 $$
 $$
 \gathered \mathcal D\big(l_1\vee\ldots\vee l_n\big)=\left[
 \op\bigcup_{i=1}^n\bigg  (\subd l_i+\sum_{j=1, j\neq i}^n\supd
 l_j\bigg  ), \ \sum_{j=1}^n \supd l_j\right],
 \\
 \mathcal D\big(l_1\wedge\ldots\wedge l_n\big) =\left[
 \sum_{j=1}^n\subd l_j, \ \op\bigcup_{i=1}^n\bigg  (\supd
 l_i+\sum_{j=1,j\neq i}^n\subd l_j\bigg) \right].
 \endgathered
 $$
For more details see \cite[Chapter~5]{SBD}):

\section{Quasidifferentiable Mappings}
\label{sec:2}

 Here we will recall the concept of quasidifferential  and consider the
 simple properties of  directional derivatives.

 \paragraph{2.1.}
 \label{par:1}
 Let $X$ be a~vector space, and let $E$ be a~$K$-space.
 Put
 $E^{\scriptscriptstyle{\bullet}}\!:=E\cup\{+\infty\}$.
 Consider  $f:X\to E^{\scriptscriptstyle{\bullet}}$
 and  $x_0\in\core(\dom(f))$.
 The record $x\in\core(C)$ means that  $C-x$ is an absorbing set.
 If, given $h\in X$, there exists
  $$
 \gathered f'(x_0)h:=
 \olim_{\alpha\downarrow0}\frac{f(x_0+\alpha h)-f(x_0)}{\alpha}\\
 =\inf_{\varepsilon>0}\sup_{0<\alpha<\varepsilon} \frac{f(x_0+\alpha
 h)-f(x_0)}{\alpha} =\sup_{\varepsilon>0}\inf_{0<\alpha<\varepsilon}
 \frac{f(x_0+\alpha h)-f(x_0)}{\alpha},
 \endgathered
 $$
 then we call  $f'(x_0)h$ the {\it one-sided
 derivative\/} or, rarely, the
 {\it Dini derivative\/} of~$f$  at
 $x_0$  in direction~$h$.
 Assume that, given $x_0$, the element  $f'(x_0)h$ exists
 for all $h\in X$. Then the mapping  $f^\prime(x_0):X\to E$ appears  which is also
 called the {\it one-sided directional derivative\/} or {\it
Dini derivative\/} at $x_0$. In this event we also say that
 $f$ is {\it directionally differentiable\/} at~$x_0$.

 \paragraph{2.2.}
 \label{par:2}
 A mapping $f$ is {\it quasidifferentiable\/}
 at~$x_0$ provided that

 (1)~$f$ is directionally differentiable at~$x_0$;

 (2)~$f'(x_0):X\to E$ is a quasilinear mapping.

 If $f$ is quasidifferentiable at~$x_0$, then the Minkowski
 duality assign to the quasilinear operator
 $f'(x_0)\in\QL(X,E)$ the element $\mathcal D \big(f'(x_0)\big)\in
 [\CSc(X,E)]$ which is called the {\it quasidifferential\/} of~$f$ at~$x_0$
 and denoted  by~$\mathcal D f(x_0)$.

 If $f'(x_0)$ can be presented as difference of sublinear operators
 $p,q\in \Sbl(X,E)$ so that $\mathcal D f(x_0)=[\partial p, \partial q]$, then
 $$
 f'(x_0)h=\sup_{S\in\partial p}S(h) -\sup_{T\in\partial q}T(h)=
 p(h)-q(h)\quad (h\in X).
 $$
 In this event $\partial p$ and $\partial q$ are respectively called
 the {\it subdifferential\/} and {\it superdifferential\/} of~$f$
 at~$x_0$ and denoted by $\subd f(x_0)$ and
 $\supd f(x_0)$. In other words,
 $$
 \mathcal Df(x_0):=[\partial p, \partial q]:=[\subd f(x_0),\supd f(x_0)].
 $$

 Assume that a quasidifferentiable mapping $f$ has   the quasidifferential at~$x_0$ of the form
  $\mathcal Df(x_0)=[\subd f(x_0),\{0\}]$
 or $\mathcal Df(x_0)=[\{0\},\supd f(x_0)]$). Then we say that $f$
 is {\it subdifferentiable\/} or, respectively,
  {\it  superdifferentiable\/}  at~$x_0$.
 If $f$ has the directional derivative $T:=f'(x_0)$ at some point $x_0\in\core(\dom(f))$
 which is a linear operator, then $f$ is subdifferentiable and superdifferentiable simultaneously
 and, moreover,  $\mathcal
 Df(x_0)=[\{T\},\{0\}]=[\{0\},\{-T\}]$.

 Each convex operator $f$ is subdifferentiable at every point
  $x_0\in\core(\dom(f))$, since we have the sublinear directional derivative
 $f'(x_0)$. In this event $\subd
 f(x_0)=\partial f(x)$. An operator $f$ is called {\it concave\/} provided that
  $-f$ is a~convex operator. A concave operator $f$
 is superdifferentiable at every point $x_0\in\core(\dom(-f))$ and, moreover,
  $\supd f(x_0)=-\partial (-f)(x_0)$. In this event the directional derivative
  $f'(x_0)$ exists too and presents a~{\it superlinear operator}; i.e.,
   $-f'(x_0)$ is a~sublinear operator.

 The differences of convex operators or, which is the same, the sums of convex and concave
  operators comprise  a much broader class of quasidifferentiable mappings.

 \paragraph{2.3.}
 Let $E$ and $F$ be  $K$-spaces.
 Consider a mapping $g:E\to F^{\scriptscriptstyle{\bullet}}$ directionally differentiable
 at~$e_0\in\core(\dom(g))$. Take $u\in E$ and $d\in F$. Suppose that
 for every sequence $(e_n)\subset E$, $e_n\downarrow 0$,
 we have
 $$
 \inf_{m\in\mathbb N}\sup_{\substack{0<\alpha<1/m\\|u'-u|\leqslant e_m}}
 \bigg|\frac{g(e_0+\alpha u')-g(e_0)}{\alpha}-d\bigg|=0.
 $$
 Then $d$ is called the {\it Hadamard derivative\/} of $g$ at~$e_0$ in direction~$u$, and
 we put $g'(e_0)u:=d$. This denotation is justified by the obvious reason that
 if the Hadamard derivative exists then so does the Dini derivative at the same point and in the same
 direction and the two derivatives coincide.
 Therefore, we can define the Hadamard derivative of~$g$ at~$e_0$ in direction
 $u$ by the formulas
 \begin{multline*}
  g'(e_0)u:=g'_{e_0}(u):=
 \inf_{m\in\mathbb N}\sup_{\substack{0<\alpha<1/m\\ |u'-u|\leqslant
 e_m}} \frac{g(e_0+\alpha u')-g(e_0)}{\alpha}
 \\
 =\sup_{m\in\mathbb N}\inf_{\substack{0<\alpha<1/m\\ |u'-u|\leqslant
 e_m}} \frac{g(e_0+\alpha u')-g(e_0)}{\alpha}.
 \end{multline*}
 If the Hadamard derivative $g'(e_0)u$ exists at $e_0$ in every direction
 $u\in E$ then we say that $g$ is {\it Hadamard differentiable\/}
at~$e_0$.

 The definition of Hadamard derivative is simplified if $F$ is a regular
 $K$-space. Recall that a~$K$-space~$F$ is called {\it regular\/}
 provided that for every nested sequence of subsets
  $F\supset  A_1\supset\ldots\supset A_n\supset\ldots$  satisfying
   $a=\inf_n\sup(A_n)$ there are finite subsets $A'_n\subset
 A_n$ enjoying the property $\olim_{n\to\infty}\sup (A'_n)=a$.

 \Proclaim{}
 Let $F$ be a~regular $K$-space. An element $d\in
 F$ is the Hadamard derivative of $g:E\to  F^{\scriptscriptstyle{\bullet}}$ at~$e_0\in\core(\dom(g))$
 in direction $u\in E$ if and only if for all sequences $(\alpha_n)\subset\mathbb R$ and $(u_n)\subset E$
 such that $\alpha_n\downarrow 0$ and
 $u_n\overset{(o)}{\longrightarrow} u$ we have
 $$
 \allowdisplaybreaks
 \gathered
 d=\olim_{n\to\infty}\frac{g(e_0+\alpha_n u_n)-g(e_0)}{\alpha_n}
 \\
 =\inf_{m\in\mathbb N}\sup_{n\geqslant m} \frac{g(e_0+\alpha_n
 u_n)-g(e_0)}{\alpha_n} =\sup_{m\in\mathbb N}\inf_{n\geqslant m}
 \frac{g(e_0+\alpha_n u_n)-g(e_0)}{\alpha_n}.
 \endgathered
 $$
 \Endproc

 \paragraph{2.4.}
 In the situation under study the Hadamard differentiability of~$g$ does not guarantee
 the continuity of the directional derivative $g'(e_0)(\cdot)$ as this happens in case
 $E=\mathbb  R^n$ and $F=\mathbb R$; cp.~\cite{34}. Let us consider the two cases
 in which the Hadamard differentiable mapping has the directional derivative continuous.
We will understand continuity of the directional derivative as follows:

 A mapping $\varphi:E\to F$ is called {\it
 $mo$-continuous\/} at~$u_0\in E$ provided that for every sequence
  $(e_n)\subset E$,  $e_n\downarrow0$, we have
 $$
 \inf_{n\in\mathbb N}\sup_{|u-u_0|\leqslant
 e_n}|\varphi(u)-\varphi(u_0)|=0.
 $$
 If $F$ is a~regular $K$-space, then $mo$-continuity
 means sequential $o$-continuity. Recall that
 $U\subset E$ is called a {\it normal\/} subset whenever
 $u_1\leqslant  e\leqslant u_2$ implies that $e\in U$ for all $u_1,u_2\in U$ and $e\in E$.

 \textbf{(1)} \proclaim{}Assume that $g:E\to
 F^{\scriptscriptstyle{\bullet}}$ is Dini differentiable at~$e_0\in\core(\dom(g))$.
 Assume further that there are a normal subset
  $U\subset E$ and an $mo$-continuous sublinear operator
 $p:E\to F$ such that $e_0\in\core(U)$ and
 $$
 |g(u_1)-g(u_2)|\leqslant p(u_1-u_2)\quad(u_1,u_2\in U).
 $$
Then $g$ is Hadamard differentiable at~$e_0$ and the directional derivative
$g'(e_0)(\cdot)$ is $mo$-continuous. \Endproc

 \textbf{(2)} \proclaim{}Assume that $F$ is a~regular $K$-space.
 If $g:E\to F^{\scriptscriptstyle{\bullet}}$
 is Hadamard differentiable at~$e_0\in\core(\dom(g))$ then
 the directional derivative $g'(e_0)(\cdot)$ is sequentially
  $o$-continuous. \Endproc

 \paragraph{2.5.}
 In~\cite{DemR2} V.~F. Demyanov and A.~M. Rubinov gave the definition of the quasidifferential
 of a mapping from a Banach space to a~Banach $K$-space. We provide a somewhat more general
 definition of quasidifferentiability.

 Let $X$ be a~topological vector space, and let $E$ be a~topological
 $K$-space. The latter means that $E$ is simultaneously
 a topological vector space and a~$K$-space having a`filter base of neighborhoods of the zero
 which consists of normal subsets of $E$. If we replace the order limit in the definition of
 directional derivative in~2.1 by the topological limit an require in the definition of quasidifferentiability
 in~2.2 that the directional derivative may be presented as difference of two continuous sublinear operators
 then we arrive at the definition of topological quasidifferential we will use later in~5.2.

We thus encounter the problem of how the two definitions of quasidifferential in 2.1 and 2.5 are related.
 Note that  if $E$ enjoys condition (A), i.e., order convergence implies topological convergence; then the two
  definitions yield the same.

 \section{Quasidifferential Calculus} \label{sec:3}

 Let us consider the main formulas for calculating the quasidifferentials of mapping with range in
 a~$K$-space. Quasidifferential calculus of  scalar functions on $\mathbb{R}^n$
 was developed by V.~F.~Demyanov, L.~N.~Polyakova, and A.~M.~Rubinov
 in \cite{DemPR} and \cite{33}; also see   \cite{32} and \cite{34}. 
 In \cite{34} there was showed how to apply the methods of quasidifferential calculus
 to the mappings with range in a Banach $K$-space. The quasidifferentials of the operators 
 with range an arbitrary $K$-space were investigated in~\cite{Bas1} and \cite{BaK1}; also see~\cite{SBD}.

 \paragraph{3.1.}
 Let  $f_1,\dots,f_n: X\to
 E^{\scriptscriptstyle{\bullet}}$ be quasidifferentiable at~$x_0\in\bigcap\nolimits_{i=1}^n\core(\dom(f_i))$. Then the sum
 $f\!:=f_1+\,\cdots+f_n$, supremum $g\!:=f_1\vee\cdots\vee f_n$, and
 infimum   $h\!:=f_1\wedge\cdots\wedge f_n$ of $f_1,\dots,f_n$ are quasidifferentiable at~$x_0$
 as well and we have the formulas:
  $$
 \gathered \mathcal D f(x_0)=\mathcal Df_1(x_0)+\,\cdots+\mathcal D
 f_n(x_0)
 \\
 =[\subd f_1(x_0)+\cdots+\subd f_n(x_0),\, \supd
 f_1(x_0)+\cdots+\supd f_n(x_0)].
 \\
 \subd g(x_0)=\bigcup_{(\alpha_1,\dots,\alpha_n)\,\in\,\Gamma_n(x_0)}
 \sum_{k=1}^n \alpha_k\Big(\subd f_k(x_0) + \sum_{l\neq k}\supd
 f_l(x_0)\Big),
 \endgathered
 $$
 $$
 \supd g(x_0)=\sum_{k=1}^n\supd f_k(x_0),\quad \subd
 h(x_0)=\sum_{k=1}^n\subd f_k(x_0),
 $$
 $$
 \gathered \supd
 h(x_0)=\bigcup_{(\alpha_1,\dots,\alpha_n)\,\in\,\Delta_n(x_0)}
 \sum_{k=1}^n \alpha_k\Big(\supd f_k(x_0) + \sum_{l\neq k}\subd
 f_l(x_0)\Big),
 \endgathered
 $$
 where
 \begin{multline*}
  \Gamma_n(x_0)\!:=\Gamma_n(x_0;f_1,\dots,f_n)\!:=
 \bigg\{(\alpha_1,\dots,\alpha_n):\ \alpha_k\in\Orth_+(E),
 \\
 \sum_{k=1}^n\alpha_k=I_E, \ \sum_{k=1}^n\alpha_k
 f_k(x_0)=f(x_0)\bigg\}.
 \end{multline*}
 \begin{multline*}
  \Delta_n(x_0)\!:=\Delta_n(x_0;f_1,\dots,f_n)\!:=
 \bigg\{(\alpha_1,\dots,\alpha_n) :\, \alpha_k\in\Orth_+(E),
 \\
 \sum_{k=1}^n\alpha_k=I_E, \ \sum_{k=1}^n\alpha_k
 f_k(x_0)=g(x_0)\bigg\}.
 \end{multline*}
  \Endproc

 \paragraph{3.2.} Let $f:X\to
 E^{\scriptscriptstyle{\bullet}}$ and $g:X\to
 \Orth(E)^{\scriptscriptstyle{\bullet}}$ be quasidifferentiable at~$x_0\in\core(\dom(f))\cap\core(\dom(g))$. 
 Then the mapping
 $gf\!=g\cdot f:X\to  E^{\scriptscriptstyle{\bullet}}$, acting by the rule $g f: x\mapsto g(x)f(x)$, 
 is quasidifferentiable at~$x_0$ too and we have the formulas
 $$
 \mathcal D(g\cdot f)(x_0)=g(x_0)\mathcal Df(x_0)+\mathcal
 Dg(x_0)f(x_0),
 $$
 with
 $$
 \aligned
  \subd(g f)(x) &= g^+(x_0)\subd f(x_0) + g^-(x_0)\supd
 f(x_0)
 \\
 &+ \subd g(x_0)f^+(x_0) + \supd g(x_0)f^-(x_0),
 \\
 \supd(g f)(x) &= g^+(x_0)\supd f(x_0) + g^-(x_0)\subd f(x_0)
 \\
 &+\supd g(x_0)f^+(x_0) + \subd g(x_0)f^-(x_0).
 \endaligned
 $$
 \Endproc

 \paragraph{3.3.} 
 We now list some conditions for the composite of quasidifferentiable mappings be
 quasidifferentiable as well.
 
 \Proclaim{Theorem.} Let $X$ be a~vector space, while $E$ and
 $F$ are $K$-spaces. Assume that $f: X\to
 E^{\scriptscriptstyle{\bullet}}$ is quasidifferentiable at~$x_0\in\core(\dom(f))$, and $g: E\to
 F^{\scriptscriptstyle{\bullet}}$ is quasidifferentiable and Hadamard differentiable
 at~$e_0:=f(x_0)\in\core(\dom(g))$ with $mo$-continuous derivative $g'(e_0)(\cdot)$. 
 Assume further that the quasidifferential $\mathcal D g(e_0)$
 is defined by the pair of order bounded
 support sets $\subd g(e_0)$ and $\supd g(e_0)$ in~$L^\sim(E,F)$. Then $g\circ
 f$ is quasidifferentiable at~$x_0$. If $\subd g(e_0)\cup\supd
 g(e_0)\subset[\Lambda_1,\Lambda_2]$ for some
 $\Lambda_1,\Lambda_2\in L^\sim(E,F)$, then
 $$
 \mathcal D(g\circ f)(x_0)= \left[\bigcup_{C\in\subd g(e_0)}\partial
 (P_C),\ \bigcup_{C\in\supd g(e_0)}\partial (P_C)\right],
 $$
 where
 $$
 P_C(x_0):=(C-\Lambda_1)\sup_{S\in\subd f(x_0)}S(x_0)
 +(\Lambda_2-C)\sup_{T\in\supd f(x_0)}T(x_0).
 $$
 \Endproc

 This formula was derived under somewhat different assumptions by~A.~M.~Rubinov; cp.~\cite{DemR2}.

  \paragraph{3.4.} Let $X$ and $Y$ be Banach spaces, while $E$ is an arbitrary
   $K$-space. A mapping $f:X\times Y\to E$
 is called {\it convex-concave\/} on~$X\times Y$ provided that
  $x\mapsto f(x,y)$ is convex at every fixed $y\in Y$,  and
 $y\mapsto f(x,y)$ is concave at every fixed $x\in
 X$. The {\it partial subdifferentials\/} $\partial_x(x_0,y_0)$ and
 $\partial_y(x_0,y_0)$ are defined as
 $$
 \aligned
 \partial_x(x_0,y_0)&\!:=\subd f_1(x_0)\!
 \\
 &=\{T\in L(X,E):\, f(x,y_0)-f(x_0,y_0)\geq T(x-x_0)\ \ (x\in X)\};
 \\
 \partial_y(x_0,y_0)&\!:=\supd f_2(y_0)\!
 \\
 &=\{S\in L(Y,E):\, f(x_0,y)-f(x_0,y_0)\leq S(y-y_0)\ \ (y\in Y)\},
 \endaligned
 $$
 where $f_1(x)\!:=f(x,y_0)$ and $f_2(y)\!:=f(x_0,y)$.

 The problem of quasidifferentiability of a convex-concave function was studied by
V.~F.~Demyanov and L.~V.~Vasilieva in~\cite{32}. In more detail, assume that $X=\mathbb{R}^n$, $Y=\mathbb{R}^m$,
and $E=\mathbb{R}$. Then a convex-concave function $f:
 \mathbb{R}^n\times\mathbb{R}^m\to \mathbb{R}^{\scriptscriptstyle \bullet}$
 is quasidifferentiable at an interior point $(x_0,y_0)$ of its domain of
 definition. Moreover,
 $$
 \mathcal D f(x_0,y_0)=[\partial_x(x_0,y_0)\times\mathbb0_m,\,
 \mathbb0_n\times\partial_y(x_0,y_0)].
 $$
 In this connection it is of interest to find some conditions for a convex-concave operator to be
 quasidifferentiable.

 \paragraph{3.5.} V.~V. Gorokhovik introduced in~\cite{Goroh1} and \cite{Goroh2}
 the concept of $\varepsilon$-quasidifferenti\-al for a real function on a finite dimensional space.
He also studies the properties of $\varepsilon$-quasidifferential; cp.~\cite{34}. 
The article~\cite{Goroh2} contains the main formulas for calculating
 $\varepsilon$-quasidifferentials and demonstrates then the necessary and sufficient
 conditions of a local extremum can be formulated for a broad class of functions in terms
 of upper and lower local approximations.
  
 This raises the problem of deriving the formulas for calculating the
  $\varepsilon$-quasidifferentials of the mapping from an arbitrary vector space to an arbitrary 
 $K$-space and studying the Gorokhovik quasidifferentiability of these mappings.
 
 \section{Necessary Optimality Conditions} \label{sec:4}

 In this section we will present the necessary conditions of an extremum for quasidifferentiable
 mappings with range in a~$K$-space, staying in the algebraic framework.
 Furthermore, we will keep the same terminology and notations as in~\cite[\S\,5.1]{SBD}. The main
 results of this section were obtained in~\cite{Bas2}; also see
 \cite[Chapter~6]{SBD}.

 \paragraph{4.1.}
 Assume that $X$ is a~vector space and
 $E$ is an arbitrary $K$-space. Consider some program
 $(C,f)$, i.e., a multiple criteria extremal problem~$x\in C,\
 f(x)\to\inf$, where $C\subset X$ is some set and $f:X\to
 E^{\scriptscriptstyle{\bullet}}$ is a mapping that is assumed quasidifferentiable
 at the appropriate point of~$\core(\dom(f))$.
 A point $x_0\in C$ is an {\it ideal local infimum (supremum)} 
 of the program $x\in C,\ f(x)\to\inf$ (or $x\in C,\ f(x)\to\sup$) provided that
 there is $U\subset X$ such that $0\in \core U$ and
 $f(x_0)=\inf\{f(x):\, x\in C\cap(x_0+U)\}$ (respectively,
 $f(x_0)=\sup\{f(x):\, x\in C\cap(x_0+U)\}$). We will understand a local extremum
 by analogy, unless specified otherwise.
 
 \paragraph{4.2.}
 Let us formulate the necessary optimality conditions in the unconstrained program, i.e., 
 in case $C=X$.

 \textbf{(1)} \proclaim{Theorem.}
 Let $f:X\to
 E^{\scriptscriptstyle{\bullet}}$ be quasidifferentiable at~$x_0\in \core(\dom(f))$. 
 If $x_0$is an ideal local optimum in the unconstrained vector program
 $f(x)\to\inf$, then $\supd
 f(x_0)\subset \subd f(x_0)$ or, which is the same, $\mathcal D
 f(x_0)\geq0$. \Endproc

 The necessary optimality conditions admit the equivalent reformulation:
  $$
 \supd f(x_0)\subset \subd f(x_0)\,\Longleftrightarrow\, 0\in
 \bigcap_{v\in \supd f(x_0)}\big(\subd f(x_0)-v\big).
 $$
 \Endproc

 \paragraph{4.3.}
 Consider some vector program of the form $(C,f)$,
 with $C:=\{x\in X:\,g(x)\leqslant 0\}$, on assuming that $f$ and $g$
are quasidifferentiable at the appropriate point. We will denote
this program by~$(g,f)$. Let us introduce the quasiregularity condition
we will use below.

 Let $X$ be an~arbitrary vector space, while $E$ and $F$ are some
 $K$-spaces. 
 Recall that a mapping $T$ from $F$ to~$E$ is a {\it Maharam operator\/} provided that
$T$ is order continuous and preserves order intervals; cp.~\cite[Chapter~4]{SBD}.
Consider $f:X\to E^{\scriptscriptstyle{\bullet}}$ and ${g:X\to
 F^{\scriptscriptstyle{\bullet}}}$. The vector program $(g,f)$
 is called {\it quasiregular\/} at~${x_0\in\core(\dom(g))}$
 provided that the following hold:

 {\bf(a)}~there are a sublinear Maharam operator  $r:F\to E$ and an
 absorbing $U\subset X$ such that $\pi_x f(x_0)\leq\pi_x f(x)$ for all $x\in
 x_0+U$, where $\pi _x:=[(r\circ g(x))^-]$ is the projection on the bang generated 
 by~$(r\circ g(x))^-$;

 {\bf(b)}~$\pi T\circ\supd
 g(x_0)\cap\pi T\circ\subd g(x_0)=\varnothing$ for all $T\in\partial r(g(x_0))$ and every nonzero
 projection $\pi\in\mathfrak P(E)$.

 Condition (a) is satisfied for instance in the case that there is
 a sublinear Maharam operator $r:F\to E$ such that  $g(x)\not\leqslant0$ implies 
 $r\circ g(x)\geq0$ for all $x\in X$.

 \paragraph{4.4.}
 \proclaim{Theorem.}
 Assume the quasiregularity condition 4.3. If a feasible point $x_0$ is an ideal
 local optimum of the quasiregular quasidifferentiable problem
 $(g,f)$, then to given $s\in\supd f(x_0)$ and $S\in\supd g(x_0)$
 there are a positive orthomorphism $\alpha\in\Orth_+(E)$ and a
 Maharam operator $\gamma\in L_+(F,E)$ such that the following system of conditions is
 compatible:
 $$
 \gathered \ker{\alpha}=\{0\},\quad \gamma\circ g(x_0)=0,
 \\
 0\in \alpha(\subd f(x_0)-s)+\gamma\circ(\subd g(x_0)-S).
 \endgathered
 $$
 \Endproc

 \paragraph{4.5.}
 Observe a few corollaries of Theorem 4.4:

 \textbf{(1)}
 \proclaim{}Assume that the quasidifferentiable vector program
  $(g,f)$ enjoys the quasiregularity condition 4.3 at a feasible
  point $x_0\in X$ and, moreover, $\pi r\circ
 g(x_0)<0$ for every nonzero projection $\pi\in\mathfrak
 P(E)$. If $x_0$ is an~ideal optimum of~$(g,f)$, then
  $ \supd f(x_0)\subset \subd f(x_0). $ \Endproc

 \textbf{(2)}~In the case of convex programs, the quasiregularity in the sense of~4.3
 agrees with the quasiregularity in the sense of~\cite[5.2.1]{SBD}. 
 Indeed, if $f$ and $\tilde g:=r\circ g$ are convex operators, then
  $$
 \subd f(x_0)=\partial f(x_0),\quad \supd f(x_0)=\{0\},\quad
 \subd\tilde g(x_0)=\partial\tilde g(x_0),\quad \supd\tilde
 g(x_0)=\{0\}.
 $$
 Furthermore,  (b) of the quasiregularity condition 4.3, i.e.,
 the equality $\pi  T\circ\supd g(x_0)\cap\pi T\circ\subd g(x_0)=\varnothing$
 valid for all $T\in\partial r$ means in this event that
 $0\notin\pi\partial\tilde g(x_0)$ which amounts to the existence of
  $h_0\in X$ satisfying $\pi\tilde  g'(x_0)h_0<0$. But
 $$
 \pi\tilde g'(x_0)h_0= \inf_{t>0}\frac{\pi\tilde
 g(x_0+th_0)-\pi\tilde g(x_0)}{t};
 $$
 hence, there are  a projection $0\neq\pi'\leq \pi$ and a real
 $t_0>0$ such that $\pi'g(x_0+t_0h_0)<\pi' g(x_0)\leq 0$. Therefore,
 the conditions of~4.3 are as follows: For some sublinear Maharam
 operator $r:F\to E$ and an absorbing set $U$, we have firstly that
 $\pi_x  f(x_0)\leq\pi_x f(x)$ for all $x\in x_0+U$ where $\pi_x:=[g(x)^-]$ 
 and, secondly, for each nonzero projection $\pi\in\mathfrak P(E)$ there are
 a projection $0\ne\pi'\leq\pi$ and $x'\in X$ with $x':=x_0+t_0h_0$
 such that $(r\circ g)(x')<0$. 

 \textbf{(3)}
 \proclaim{}A feasible point
 $x_0\in\core(\dom(f))\cap\core(\dom(g))$ is an ideal optimum
 of a~quasiregular convex program $(g,f)$ if and only if
 there are a positive orthomorphism
 $\alpha\in\Orth_+(E)$ and a Maharam operator $\gamma\in L_+(F,E)$ such that
 the following system of conditions is compatible:
 $$
 \ker{\alpha}=\{0\},\quad \gamma\circ g(x_0)=0,\quad
 0\in\alpha\partial f(x_0)+\gamma\circ\partial g(x_0).
 $$
 \Endproc

 \textbf{(4)}~Let $f,\varphi:X\to
 E^{\scriptscriptstyle{\bullet}}$ and $g,\psi:X\to
 F^{\scriptscriptstyle{\bullet}}$ be quasidifferentiable at an appropriate point.
 Reduce each of the extremal problems
 $$
 \gathered \psi (x)\geq 0,\quad f(x)\to \inf;
 \\
 g(x)\leq 0,\quad \varphi (x)\to \sup;
 \\
 \psi (x)\geq 0,\quad \varphi (x)\to \sup
 \endgathered
 $$
 to the above problem $(g,f)$ by letting $g:=-\psi$ and $f:=-\varphi$. 
 In this event we encounter some obvious modifications of the quasiregularity
 condition. Indeed, the quasiregularity condition for the program
 $\psi (x)\geq  0,\ f(x)\to \inf$ means the existence of a sublinear Maharam
 operator $r:F\to E$ such that, firstly,,
 $\pi_x f(x_0)\leq\pi_x f(x)$ for all $x\in X$ with $\pi _x:=[(r\circ\psi(x))^+]$ 
 and, secondly,  $\pi T\circ\subd\psi(x_0)\cap\pi
 T\circ\supd\psi(x_0)=\varnothing$ for every
 $T\in\partial r(\psi (x_0))$ and every nonzero projection $\pi\in\mathfrak
 P(E)$.

 \Proclaim{}If a feasible point $x_0$ is an ideal local optimum of
 a general quasiregular quasidifferentiable program
 $\psi (x)\geq 0,\ f(x)\to \inf$, then to all $s\in\supd f(x_0)$ and
 $S\in\subd\psi(x_0)$ there are a positive orthomorphism
  $\alpha\in\Orth_+(E)$ and a Maharam operator $\gamma\in L_+(F,E)$ such that the
  following system of conditions is compatible:
 $$
  \gathered \ker{\alpha}=\{0\},\quad
 \gamma\circ\psi(x_0)=0,
 \\
 0\in \alpha(\subd f(x_0)-s)+\gamma\circ(\supd \psi(x_0)-S).
 \endgathered
  $$
 \Endproc

 {\bf(5)}~Theorem 4.4, the main result of 2.4, was obtained by
 E.~K.~Basaeva in~\cite{Bas2}; also see~\cite{SBD}.
 In the scalar case $E=F=\mathbb R$ and $X=\mathbb R^n$ Theorem~4.4
 is well known; for instance, see \cite[Theorem V.3.2]{SBD}.
 In this event the quasiregularity condition may be slackened, but this will lead to
 slackening the necessary optimality conditions. In more detail, if the
 closure of $\{h\in X:\, g'(x_0)h<0\}$ is $\{h\in X:\, g'(x_0)h\leqslant 0\}$ (which is regularity), 
 then
 $$
 0\in(\subd f(x_0)-s)+\cl\cone(\subd g(x_0)-S),
 $$
 for all
 $s\in\supd f(x_0)$ and $S\in\supd g(x_0)$, where
  $\cl\cone(\mathcal U)$ stands for the closed cone hull of~$\mathcal U$.

  \paragraph{4.6.}
 The quasiregularity condition 4.3 allows us to write down the necessary
 optimality conditions for $S\in\supd g(x_0)$ provided that
 $\pi T\circ\supd g(x_0)\cap\pi
 T\circ\subd g(x_0)=\varnothing$ for every $T\in\partial r(g(x_0))$ and every nonzero
 $\pi\in\mathfrak P(E)$. If the last condition is not fulfilled then there is
 a maximal projection $\rho$ such that the condition is valid for all
 projections $\pi\leq \rho$, and part of the necessary conditions may be written as 
 in~4.4 but with a~constraint on~ $\rho$. Namely,
 $$
 0\in \rho \alpha(\subd f(x_0)-s)+\rho\gamma(\subd g(x_0)-S).
 $$
 The necessary optimality  conditions at~$\rho^d$ will principally different:
 To all $s\in\supd f(x_0)$ and $S\in\supd  g(x_0)$ there is a positive orthomorphism
  $\alpha\in\Orth_+(\rho^d E)$, $\ker{\alpha}=\{0\}$ such that
 $$
 0\in \rho^d\alpha(\subd f(x_0)-s)+ \rho^d\Cop(\partial
 r(g(x_0))\circ\subd g(x_0)-S),
 $$
 were $\Cop(\mathcal U)$ stands for $\cl(\mix(\co\mathcal
 U))$, while $\co(A)$ is the convex hull of $A$, while $\mix(A)$ is the collection of all mixings
 of the elements of $A$ by all partitions of unity in
  $\mathfrak{P}(E)$ (cp.~\cite[Appendix~4]{SBD}), and $\cl(A)$ is the closure of $A$ with respect to pointwise 
  $o$-convergence. This fact can be established on using the Vector Bipolar Theorem
   \cite[p.~190]{SBD}; also see~\cite{KVD}.

 \section{Accounting for the Containment Constraints}
 \label{sec:5}

Let us consider the necessary optimality conditions in the case that
the program under study involves the constraint that the solution must belong to
a given set. The regularity condition on the set is convenient to formulate in terms of
the topology of the ambient vector space. In other words, we have to define
topological quasidifferentials. To this end it suffices to change the scope of the
concept of quasilinear mapping so as to imply now that
 a quasilinear operator must admit representation as difference of 
 {\it continuous\/} sublinear operators.

 \paragraph{5.1.}
 Let $X$ be a~topological vector space, while
 $E$ is a~topological $K$-space an $A^c$ is the algebra of
 continuous orthomorphisms of~$E$. Assume that the positive cone
 of a~ topological $K$-space is normal by default.
 In this event the Minkowski duality $\partial$ defines a bijection
 between the sets of  (total = everywhere defined) continuous sublinear
operators and the collection of equicontinuous support sets; cp.~\cite[3.2.2\,(1)]{SBD}.

 Let $\tQL(X,E)$ stand for the part of $\QL(X,E)$ consisting of the
 quasilinear operators that are presentable as differences of continuous sublinear operators.
 Clearly, $\tQL(X,E)$ is a~lattice ordered $A^c$-module. The module and lattice
 operations as well as order, are induced from $\QL(X,E)$. We call the members of $\tQL(X,E)$
 {\it continuous quasilinear operators}.

 By analogy, the collection of equicontinuous support sets
 $\tCSc(X,E)$ is defined as the part of $\CSc(X,E)$ consisting of
 the support sets of sublinear operators; cp.~\cite[.2.2\,(1)]{SBD}. 
The relevant restriction of the isomorphism $\mathcal D$ of~1.13 will be
denoted by~$\mathcal D^c$. Clearly, $\mathcal
 D^c$ is an isomorphism between the $A^c$-modules $[\tCSc(X,E)]$ and
 $\tQL(X,E)$.

 \paragraph{5.2.}
Therefore, to preserve the formulas of quasidifferential calculus of
\S\,4 in the topological setting it suffices ti require
that the directional derivative in the definition of~2.2 may be presented
as difference of continuous sublinear operators.

Let $X$ be a~topological vector space and let $E$ be a~topological
$K$-space. Consider $f:X\to
 E^{\scriptscriptstyle{\bullet}}$ and $x_0\in\core(\dom(f))$.
 We say that $f$ is {\it topologically quasidifferentiable\/} at~$x_0$ provided that
 the directional derivative $f'(x_0)$ exists  at~$x_0$ and presents
 a continuous quasilinear operator.
 
 Therefore, if $f$ is topologically quasidifferentiable at~$x_0$, then
 the Minkowski duality put in correspondence to the quasilinear operator $f'(x_0)\in\tQL(X,E)$ 
 the element $\mathcal D  (f'(x_0))\in [\tCSc(X,E)]$ called the {\it  topological
 quasidifferential\/} of~$f$ at~$x_0$ and denoted by~$\mathcal D^cf(x_0):=[\tsubd f(x_0),\tsupd f(x_0)]$. 
 Here
 $\tsubd f(x_0)$ and $\tsupd f(x_0)$ stand respectively for the {\it
topological subdifferential\/} and {\it topological superdifferential\/} of~$f$ at~$x_0$.

 The formulas comprising the calculus of topological quasidifferentials
 coincide with their algebraic analogs in~2.2 if we replace  $\mathcal D$
 with~$\mathcal D^c$.

 \paragraph{5.3.}
 Let $X$ be a~topological vector space, while
 ${C\subset X}$ and $x_0\in C$. The {\it feasible direction cone\/}
 $\Fd(C,x_0)$ of~$C$ at~$x_0$ is introduced as follows:
 $$
 \Fd(C,x_0):=\{h\in X:\ (\exists\varepsilon>0)\,
 x_0+[0,\varepsilon)h\subset C\}.
 $$
 We say that $C$ is {\it $K$-regular\/} at~$x_0$ whenever
 $K$ is a convex cone and $K\subset\cl\big(\Fd(C,x_0)\big)$. Given
 some set $C$ that is $K$-regular at~$x_0$, we introduce a {\it normal cone\/} of~$C$ at~$x_0$
 as follows: $N_E(C,x_0):=\pi_E(K):= \{T:\, Tk\leq0,\, k\in K\}$. Obviously,
 this definition yields a nonunique cone.
 
 \Proclaim{}If $C\subset X$ is $K$-regular at~$x_0\in C$ and $f:X\to E^{\scriptscriptstyle{\bullet}}$
 is quasidifferential at $x_0$, where $x_0\in \core(\dom(f))$. For
 $x_0$ to be an ideal local optimum of~$(C,f)$ it is necessary that
  $$
 \tsupd f(x_0)\subset\tsubd f(x_0)+ N_E(C,x_0).
 $$
 \Endproc

 \paragraph{5.4.}
 Consider a vector program $(g,f)$ with the extra constraint $x\in C$ for some
 $C\subset X$. We denote this program by $(C,g,f)$. Assume that $x_0\in C\cap\core(\dom(f))\cap\core(\dom(g))$ and suppose that
 $f:X\to E^{\scriptscriptstyle{\bullet}}$ and $g:X\to
 F^{\scriptscriptstyle{\bullet}}$ are topologically quasidifferentiable at~$x_0$.
 The vector program $(C,g,f)$ is called {\it quasiregular\/} at~$x_0$ provided that

 {\bf(a)}~there are a continuous Maharam operator $r:F\to
 E$ and a neighborhood $U$ of~$x_0$ such that $\pi_x f(x_0)\leq\pi_x f(x)$
 for all $x\in C\cap  U$, where $\pi_x:=[(r\circ
 g(x))^-]$ is the projection to the band generated by~$(r\circ g(x))^-$;

 {\bf(b)}~$C$ is $K$-regular at~$x_0$;

 {\bf(c)}~for every $T\in\partial r(g(x_0))$ and every nonzero
 $\pi\in\mathfrak P(E)$ we have
 $$
 \pi
 T\circ\tsupd g(x_0)\cap\big(\pi T\tsubd g(x_0)+ \pi
 N_E(C,x_0)\big)=\varnothing.
 $$

 \paragraph{5.5.}
 \proclaim{Theorem.}
 Let $f$ and $g$ be quasidifferentiable at~$x_0\in
 C\cap\core(\dom(f))\cap\core(\dom(g))$, and let the vector program
 $(C,g,f)$ be quasiregular at~$x_0$. If $x_0$ is an~ideal local optimum
 of~$(C,g,f)$, then to all $s\in\tsupd
 f(x_0)$ and $S\in\tsupd g(x_0)$ there are a continuous orthomorphism
 $\alpha\in\Orth(E)$, a continuous Maharam operator $\gamma\in
 L_+(F,E)$, and a continuous linear operator $\lambda\in L(X,E)$ such that
 the following system of conditions is compatible:
 $$
 \gathered 0\leq\alpha\leq I_E,\quad \ker{\alpha}=\{0\},\quad
 \lambda\in N_E(C,x_0),\quad \gamma\circ g(x_0)=0,
 \\
 -\lambda\in \alpha(\tsubd f(x_0)-s)+\gamma\circ(\tsubd g(x_0)-S).
 \endgathered
 $$
 \Endproc

 \paragraph{5.6.}
 We call  $\{x_1^0,\dots, x_n^0\}\subset C$ a~{\it generalized local optimum\/} 
 of a program~$(C,f)$ provided that there is a neighborhood $U$ of the zero
 such that $f(x_1^0)\wedge
 \ldots\wedge f(x_n^0)\leq f(x_1)\wedge \ldots\wedge f(x_n)$ for all
 $x_i\in (x_i^0+U)\cap C$ and $i:=1,\dots,n$.

 \Proclaim{}Let $f:X\to
 E^{\scriptscriptstyle{\bullet}}$ be quasidifferentiable at each of the
 feasible points $x_1^0,\dots, x_n^0\in\core(\dom(f))$. If
 $\{x_1^0,\dots, x_n^0\}$ is a generalized local optimum of the unconstrained program
  $f(x)\to\inf$, then for all  $\alpha_1,\dots,\alpha_n\in\Orth_+(E)$ such that
 $$
 \alpha_1+\,\ldots+\alpha_n=I_E, \quad \sum_{i=1}^n\alpha_i
 f(x_i^0)=f(x_1^0)\wedge\ldots\wedge f(x_n^0),
 $$
 we have the inclusions
 $$
 \alpha_k \tsupd f(x_k)\subset \alpha_k \tsubd f(x_k)\quad
 (k:=1,\dots,n).
 $$
 \Endproc

 \paragraph{5.7.}
 \proclaim{Theorem.}
 Let $f:X\to
 E^{\scriptscriptstyle{\bullet}}$ be quasidifferentiable at~$x_1^0,\dots,x_n^0\in\core(\dom(f))$
 and let $C\subset X$ be  $K_l$-regular at~$x_l^0\in C$ for all $l:=1,\dots,n$, where
  $K_1,\dots,K_n$ are convex cones. If $\{x_1^0,\dots,
 x_n^0\}$ is a generalized local optimum of the program
  $(C,f)$, the for all $\alpha_1,\dots,\alpha_n\in\Orth_+(E)$ such that
   $$
 \alpha_1+\,\ldots+\alpha_n=I_E, \quad \sum_{k=1}^n\alpha_k
 f(x_k^0)=f(x_1^0)\wedge\ldots\wedge f(x_n^0),
 $$
 we have the inclusions
 $$
 \alpha_k \tsupd f(x_k^0)\subset \alpha_k \tsubd
 f(x_k^0)+N_E(C,x_k^0) \quad (k:=1,\dots,n).
 $$
 \Endproc

 \paragraph{5.8.}
 Consider some vector program $(C,g,f)$. Assume that
 $$
 x_i^0\in C\cap\core(\dom(f))\cap\core(\dom(g)),
 $$
  while
  $f:X\to E^{\scriptscriptstyle{\bullet}}$ and $g:X\to
 F^{\scriptscriptstyle{\bullet}}$ are topologically
 quasidifferentiable at~$x_i^0$ for all $i:=1,\dots,n$. The vector program
  $(C,g,f)$ is called {\it quasiregular\/} at~$\{x_1^0,\dots,x_n^0\}$
  provided that the following conditions are fulfilled:

 {\bf(a)}~there are a continuous Maharam operator $r:F\to
 E$ and neighborhoods $U_i$ of~$x_i^0$ such that 
 $\pi_x e\leq\pi_x f(x)$ for all $x\in  C\cap U_i$, where
 $e:=f(x_1^0)\wedge\ldots\wedge f(x_n^0)$ and $\pi_x:=[(r\circ
 g(x))^-]$ is the projection on the band generated by~$(r\circ  g(x))^-$;

 {\bf(b)}~$C$ is $K_i$-regular at~$x_i^0$;

 {\bf(c)}~$\pi  T\circ\tsupd g(x_i^0)\cap\big(\pi T\circ \tsubd g(x_i^0)+ \pi
 N_E(C,x_i^0)\big)=\varnothing$ for all $T\in\partial r(g(x_i^0))$, $i:=1,\dots,n$, and $\pi\in\mathfrak
 P(E)$.

 \paragraph{5.9.}
 \proclaim{Theorem.}
 Assume that $g:X\to F^{\scriptscriptstyle{\bullet}}$ and ${f:X\to
 E^{\scriptscriptstyle{\bullet}}}$ are quasidifferentiable at~$x_1^0,\dots, x_n^0\in C\ \cap\ \core(\dom(f))\ \cap\ \core(\dom(g))$. Assume further that the vector program
 $(C,g,f)$ is quasiregular in the sense of~$5.8$ at~$\{x_1^0,\dots, x_n^0\}$. 
 If $\{x_1^0,\dots, x_n^0\}$ is a generalized local optimum of~$(C,g,f)$, then
 to all $s_i\in \supd f(x_i^0)$ and $S_i\in\supd g(x_i^0)$ there are orthomorphisms
  $\alpha_1,\dots,\alpha_n\in\Orth(E)$, continuous Maharam operators $\gamma_1,\dots,\gamma_n\in L_+(F,E)$,
   and continuous linear operators $\lambda_i\in L(X,E)$ such that
 $$
 \gathered 0\leqslant\alpha_i\leqslant I_E,\quad
 \ker(\alpha_1)\cap\ldots\cap\ker(\alpha_n)=\{0\},
 \\
 \gamma_i\circ g(x_0)=0,\quad\lambda_i\in N_E(K_{\xi_i}),
 \\
 -\lambda_i\in \alpha_i(\tsubd f(x_i^0)-s_i) +\gamma_i\circ(\tsubd
 g(x_i^0)-S_i)\quad (i:=1,\dots, n).
 \endgathered
 $$
 \Endproc

 \paragraph{5.10.} 
 It seems worthwhile to study necessary optimality conditions
 for multiple criteria extremal problems of the type $(C,g,f)$ with
 $g$ and $f$ integral quasidifferentiable operators




\noindent
{\it Elena K.~Basaeva}
\par\smallskip
{\leftskip\parindent\small
\noindent
Southern Mathematical Institute\\
 22 Markus Street\\
Vladikavkaz, 362027, RUSSIA\\
E-mail: helen@smath.ru
\par}

\noindent
{\it Anatoly G.~Kusraev}
\par\smallskip
{\leftskip\parindent\small
\noindent
Southern Mathematical Institute\\
 22 Markus Street\\
Vladikavkaz, 362027, RUSSIA\\
E-mail: kusraev@smath.ru
\par}

\noindent
{\it Sem\"en S.~Kutateladze}
\par\smallskip
{\leftskip\parindent\small
\noindent
Sobolev Institute of Mathematics\\
4 Koptyug Avenue\\
Novosibirsk, 630090, RUSSIA\\
E-mail: sskut@math.nsc.ru
\par}

\end{document}